\documentclass[9pt]{article}
\newtheorem{theorem}{Theorem}
\newtheorem{question}{Question}
\def\prob{\mbox{Prob}}
\def\expect{{\mathbf E}}
\def\grad{\nabla}
\def\unif{{\overline{p}}}
\def\qed{{\rm QED}}

\begin{document}

\vspace{2in}
\begin{center}
{\large \bf A Dynamic Model of Social Network Formation} \\
\end{center}
\vspace{5ex}
\begin{flushright}
Brian Skyrms \footnote{ School of Social Sciences, University of
California at Irvine, Irvine, CA 92607} \\

Robin Pemantle \footnote{Research supported in part by National
Science Foundation grant \# DMS 9803249}$^,$\footnote{Department
of Mathematics, The Ohio State University, 231 W. 18 Ave.,
Columbus, OH 43210}

\end{flushright}

{\bf ABSTRACT:} \break
 We consider a dynamic social network model in which agents play
repeated games in pairings determined by a stochastically evolving
social network.  Individual agents begin to interact at random,
with the interactions modeled as games.  The game payoffs
determine which interactions are reinforced, and the network
structure emerges as a consequence of the dynamics of the agents'
learning behavior.  We study this in a variety of game-theoretic
conditions and show that the behavior is complex and sometimes
dissimilar to behavior in the absence of structural dynamics.  We
argue that modeling network structure as dynamic increases realism
without rendering the problem of analysis intractable.

\setcounter{equation}{0}
\section{Introduction}

Pairs from among a population of ten individuals interact
repeatedly.  Perhaps they are cooperating to hunt stags and
rabbits, or coordinating on which concert to attend together;
perhaps they are involved in the somewhat more antagonistic
situation of bargaining to split a fixed payoff, or attempting to
escape the undesirable but compelling equilibrium of a Prisoner's
Dilemma. As time progresses, the players adapt their strategies,
perhaps incorporating randomness in their decision rules, to suit
their environment. But they may also exert control over their
environment.  The players may have choice over the pairings,
though not perfect information about the other players.  They may
improve their lot in two different ways.  A child who is being
bullied either learns to fight better, or to run away. Similarly,
a player who obtains unsatisfactory results may choose either to
change strategies or to change associates. Regardless of whether
the interactions are mostly cooperative or mostly antagonistic, it
is natural and desirable to allow evolution of the social network
(the propensity for each pair to interact) as well as the
individuals' strategies.

We build a model that incorporates both of these modes of
evolution. The idea is simple.
\begin{quote} (*) \hfill \\
Individual agents begin to interact at random.  The interactions
are modeled as games.  The game payoffs determine which
interactions are reinforced, and the social network structure
emerges as a consequence of the dynamics of the agents' learning
behavior.
\end{quote}
As the details of the specific game and the reinforcement dynamics
vary, we then obtain a class of models.  In this paper we treat
some simple reinforcement dynamics, which may serve as a base for
future investigation.

The idea of simultaneous evolution of strategy and social network
appears to be almost completely unexplored.  Indeed, the most
thoroughly studied models of evolutionary game theory assume {\em
mean-field} interactions, where each individual is always equally
likely to interact with each other.  Standard treatments of
evolutionary game dynamics (1,2) operate entirely in this
paradigm. This is due, to a large extent, to considerations of
theoretical tractability of the model.  Models have been
introduced that allow the agents some control over their choice of
partner (3), but the control is still exerted in a mean-field
setting: one chooses between the present partner and a new pick at
random from the whole population.

Evolutionary biologists know that evolutionary dynamics can be
affected by non-random encounters or population structure, as in
Sewall Wright's models of assortative mating (4). Wright (5)
already realized that positive correlation of encounters could
provide an account of evolution of altruism. Thus the need for
social network models has been long recognized.

When the social network is modeled, it is almost always
static\footnote{An exception, perhaps, is a preprint we have
recently encountered  by Jackson and Watts (6)}.  Interactions,
for example, may be posited to occur only between players whose
locations are close, according to some given spatial data.
Biological models in which encounters are governed by spatial
structure have become increasingly frequent in the 1990's; see for
example the work of Durrett, Levin and Neuhauser (7,8,9).  A
similar hypothesis of spatial structure, in a game theory context,
arises in (10). Here, technology from statistical mechanics is
adapted to the analysis of games whose interactions take place
between neighbors in a grid.

A number of recent investigations by game theorists, some directly
inspired by biological models, have shown that the dynamics of
strategic interaction can be strikingly different if interaction
is governed by some spatial structure, or more generally, some
graph structure (11,12,13).  For instance, one-shot Prisoner's
Dilemma games played with neighbors on a circle or torus allows
cooperation to evolve in a way that the random encounter model
does not.  The spatial or graph structure can be important in
determining which equilibria are possible, whether repeated
interactions can be expected to converge to equilibrium, and if
so, how quickly convergence takes place (14).

Since the outcome of a repeated game may vary with the choice of
network model, it is important to get the network model right.
Further progress in the theory of games and adaptive strategies
would be greatly enhanced by a theory of networks of social
interaction. In particular, it would be desirable to have a
framework within which models may be developed that are both
tractable, and plausible as a mechanism governing interactions
among a population of agents seeking to improve their lot.

When the network changes much more slowly than do the strategies
of individuals, it is reasonable to model the social network by a
structure which is fixed, though possibly random. The question of
realistically modeling the randomness in such a case is taken up
in a number of papers, of which a recent and well known example is
the ``small world'' model (15).  In the other extreme (16,17,18)
evolution of social structure is modeled by agents moving on a
fixed graph in the absence of strategy dynamics.

In the general case, however, interaction structures are fluid and
evolve in tandem with strategy.  What is required here is a
dynamics of interaction structure to model how social networks are
formed and modified. We distinguish this {\em structure dynamics}
from the {\em strategic dynamics} by which individuals change
their individual behaviors or strategies.

In this paper we introduce a simple, additive model for structure
dynamics, and explore the resulting system under several
conditions: with or without discounting of the past, with or
without added noise, and in the presence or absence of strategic
dynamics.  Common to all our models is a stochastic evolution from
a (usually symmetric) initial state.  Individuals in a population
start out choosing whom to interact with at random, then modify
their choices according to how their choice is reinforced, and the
process is repeated.  An infinite variety of such models is
possible.  We will consider only a few basic models, meant to
illustrate that rigorous results on structure dynamics are not out
of reach, and that further inquiry will be profitable.

We first consider a baseline case of uniform reinforcement.  Here,
any choice of partner is reinforced as strongly as any alternative
choice would have been.  In other words, the interaction game
between any pair of players always produces a constant reward or
punishment.  One might expect that such cases would not lead to
interesting dynamics, but that is far from the truth.  We show
both by computer simulation and analytically how structure emerges
spontaneously even in these cases.  Since the strategic dynamics
here are trivial, the baseline case is intended mostly as a
building block on which more interesting strategic dynamics are to
be grafted.  We note, however, that the constant reward game is
not completely unreasonable.  Studies have shown that in the
absence of other environmental attributes, sheer familiarity
brings about positive attitudinal change (19).  In fact, an
abstract model of network evolution under uniform positive
re-weighting has appeared before under the name of ``Reinforced
Random Walk'' (20).

Next, we move to the case where players of different types play a
non-trivial game and are reinforced by the payoffs of the game.
Here, we examine the co-evolution of behavior and structure when
the structural dynamics and strategic dynamics are both operative.
The relative speeds of structural dynamics and strategic dynamics
affect which equilibrium is selected in the game.  In particular,
this can determine whether the risk-dominant or payoff-dominant
equilibrium is selected.

\section{Making friends: a baseline model of uniform reinforcement}

\subsection{Friends I: asymmetric weights}

Each morning, each agent goes out to visit some other agent.  The
choice of whom to visit is made by chance, with the chances being
determined by the relative {\em weights} each agent has assigned
to the others. For this purpose, agent number $i$ has a vector of
weights $\langle w_{i1} , \ldots , w_{in} \rangle$ that she
assigns to other players (assume $w_{ii} = 0$).  Then she visits
agent $j$ with probability
\begin{equation} \label{eq:prob1}
\prob (\mbox{agent } i \mbox{ visits } j) = {w_{ij}
   \over \sum_k w_{ik}} \, .
\end{equation}
Here we are interested in a symmetric baseline model, so we will
assume that all initial weights are 1.  Initially, for all agents,
all possible visits are equiprobable.

Every agent is treated nicely on her visit and all are treated
equally nicely.  They each get a reinforcement of 1.  Each agent
then updates her weight vector by adding 1 to the weight
associated with the agent that she visited.  Her probabilities for
the next round of visits are modified accordingly.  At each stage
we have a matrix $p_{ij}$ of probabilities for $i$ to visit $j$.
Do these probabilities converge, and if so to what?

Given all the symmetry built into the starting point and the
reinforcement, it is perhaps surprising that all sorts of
structure emerge.  Here is a description of a simulated sample run
of length 1,000.  The probabilities, to two decimal places, seem
to converge after a few hundred rounds of visits, to a matrix that
is anything but uniform (and to a different matrix each time the
process is run from the initial, symmetric weights).  There is one
agent, A, who visits another agent, B, more than half the time.
There is no reciprocation, so this has no bearing on how often B
visits A, and in fact most agents will not visit any one agent
more than a third of the time.

In the analysis section, we show that this outcome is typical.
\begin{theorem} \label{theorem:friends1}
The probability matrix for Friends I with $n$ players will
converge to a random limit $p$ as time goes to infinity.  The
distribution of the limit is that the rows of $p$ are independent,
each having Dirichlet distribution (ignoring the zero entry on the
diagonal) whose parameters are $n-1$ ones.
\end{theorem}
Thus we see spontaneous emergence of structure.  This type of
simple model has been used before in the economics literature to
explain the stabilization of market shares at seemingly random
equilibria, due to random reinforcement in the early phases of
growth of an industry (21). We remark that the choices made by
each agent are independent of the choices made by each other
agent, so the social aspect of the model is somewhat degenerate
and the model may be viewed as a model of individual choice.
Nevertheless, it fits our definition of social network model in
that it gives a probabilistic structure to interactions; one may
then extend the model so the interactions are nontrivial games.

\subsection{Friends II: symmetrized reinforcement}

Suppose now that the interaction is as pleasant to the host as the
visitor.  Thus when agent $i$ visits agent $j$, we add 1 to both
$w_{ij}$ and $w_{ji}$.  A typical outcome for 10 agents after
1,000 rounds of visits looks similar to the table for Friends I,
except that the entries are nearly symmetric.  There are, however,
subtle differences that may cause the two models to act very
differently when strategic dynamics are introduced.  To see these
differences, we describe what is typically observed after 10 runs
of a simulation of Friends II to time 1,000 for a set of three
agents, this being the minimum population size for which
structural dynamics are interesting.  What we see typically is one
or two runs in which each players visits are split evenly (to two
decimal places) between the others.  We see another several runs
that are close to this.  We see one run or so in which two agents
nearly always visit the third agent, which splits its time among
the other two.  The remaining runs give something between these
extreme outcomes.

What may not be apparent from such data is that the limiting
weights for Friends II are always $1/2$.  Only a small fraction of
sample outcomes decisively exhibit the proven limiting behavior.
The data, in other words, show that after 1,000 iterations, the
weights may still be far from their limiting values; when this is
the case, one of the three agents is largely ignored by the other
two, and visits each of the other two herself equally often. Since
the lifetime of many adaptive games is 1,000 rounds or fewer, we
see that limiting behavior may not be a good guide to behavior of
the system on time scales we are interested in.  The analysis
section discusses both limiting results for this model and finite
time behavior.  When the population size is more than 3, the
weights will always converge, but the limit is random and
restricted to the subspace of symmetric matrices.  Again,
convergence of the weights to their limiting values is slower than
in the non-reciprocal game of Friends I.
\begin{theorem} \label{theorem:friends2}
The probability matrix $p_{ij}$ for Friends II with $n$ players
converges to random limit $p$ as time goes to infinity.  If $n =
3$, the limit is the matrix all of whose off-diagonal entries are
$1/2$.  In general, the limit may be any symmetric matrix whose rows
sum to 1; that is, the closed support of the random limit is the
entire subspace of symmetric stochastic matrices.
\end{theorem}

\subsection{Analysis of Friends I and II}

To fit this in the framework of~$(*)$, construct the following
degenerate games.  Each of the two players has only one strategy,
and the payoff matrix is as follows.

\hspace{.5in}
\begin{tabular}{c|c|}
Friends I & Host \\
\hline
Visitor & $(1,0)$ \\
\hline
\end{tabular}
\hspace{.5in}
\begin{tabular}{c|c|}
Friends II & Host \\
\hline
Visitor & $(1,1)$ \\
\hline
\end{tabular}

The weights $w_{ij}$ are initialized to 1 for $i \neq j$, and are
then updated according to
\begin{equation} \label{eq:update weights}
w_{ij} (t+1) = w_{ij} (t) + u(i,j;t)
\end{equation}
where $w_{ij} (t)$ is the weight agent $i$ gives to agent $j$ at
time $t$ and $u(i,j;t)$ is the payoff of the game played at time
$t$ between visitor $i$ and host $j$ (and zero if this visit did
not occur at time $t$).  This, together with specification of the
visitation probabilities in equation~(\ref{eq:prob1}), defines the
model. Changing the initial weights does not affect the
qualitative behavior of any model, so there is no need to vary the
initialization.

For Friends I, the updating of the weights for any one agent is
the same as a P\'olya urn process (22). Each agent can be thought
of as having an urn with balls of $n-1$ colors, one color
representing each other agent.  Initially there is one ball of
each color in the urn.  The agent picks a ball at random,
indicating whom she should visit, then returns it to the urn along
with an extra ball of the same color.  The urns belonging to
different agents are statistically independent.

The analysis of this process is well known (23, Chapter 4). It is
easy to show that the sequence of draws for each agent is {\em
exchangeable}, that is, permuting a sequence does not change its
probability.  Hence by the de Finetti representation theorem, the
random sequence of draws from an urn is equivalent to a mixture of
multinomial processes, that is, of sequences of independent draws.
The mixing measure is easily seen to be Dirichlet.  Consequently,
the visiting probabilities converge with probability one, but they
can converge to anything. That they converge to the uniform
vector, where each agent has equal probability to visit each
other, has prior probability zero.

Furthermore, convergence to the limiting probability matrix is
quite rapid.  Let $p (t)$ denote the matrix whose $(i,j)$-entry is
$p_{ij} (t)$. Then exchangeability implies that, conditional on
the limit matrix $p = \lim_{t \to \infty} p (t)$, the sequence of
visits is a sequence of independent, identically distributed draws
from the limit distribution.  Thus at time $t$, the central limit
theorem implies that $p (t) - p$ is $t^{-1/2}$ times a
multivariate normal.

For Friends II, exchangeability fails.  This is not surprising,
since the property of exchangeability is not very robust.  More
surprising, however, is that the sequence of probability matrices
$p(t)$ does not form a martingale.  To explain this terminology,
let $\expect_t$ denote the expectation conditioned on the values
at time $t$.  A simple computation shows that for Friends I, the
expected value of $p_{ij} (t+1)$ conditioned on the time $t$ value
is equal to $p_{ij} (t)$: since $w_{ij}$ increases only when $i$
visits $j$, we have
\begin{eqnarray*}
\expect_t p_{ij} (t+1) & = & \expect_t \sum_{k=1}^n p_{ik} (t) {w_{ij} +
   \delta_{jk} \over 1 + \sum_{l=1}^n w_{il} (t)} \\
& = & {w_{ij} (t) + p_{ij} (t) \over 1 + \sum_{l=1}^n w_{il}(t)} \\
& = & p_{ij} (t) \, .
\end{eqnarray*}
Even without exchangeability, the martingale convergence theorem
(24, Section~4.2) implies convergence of the quantities
$p_{ij}$, though it says very little about the limit.

For Friends II the complete analysis may be found in (25).  Here
is an outline of what is found there.  A computation similar to
the one for Friends I shows that
$$\expect_t p (t+1) = p (t) + {1 \over t} F(p(t))$$
where $F$ is a certain function on symmetric $n$ by $n$ matrices.
In other words, the random sequence of matrices $\{ p (t) : t = 1
, 2 , \ldots \}$ is a {\em stochastic approximation} in the sense
of Robbins and Monro (26), driven by the vector field
$F$.  General results of (27) and (28) now imply that
$p (t)$ converges to the set where $F$
vanishes.  To show that $p (t)$ always converges to a single
point, Pemantle and Skyrms (25) compute a Lyapunov function for
$F$, that is, a function $V$ for which $\grad V \cdot F < 0$ with
equality only when $F = 0$.  This, together with an efficiency
inequality (bounding the angle between $f$ and $\grad V$ away from
ninety degrees), establish convergence of $p$.  The remainder of
Theorem~\ref{theorem:friends2} is then established by showing the
only stable zeros of the vector field $F$ are the symmetric
matrices with row sums all equal to 1, and that the possible limit
points of $p (t)$ are exactly the stable equilibria of the flow
determined by $F$.

Determination of the rate of convergence of $p(t)$ to its limit is
somewhat different in this case.  Due to the presence of unstable
equilibria from the flow determined by $F$, there is a possibility
of being stuck near one of these equilibria for a long time before
eventually following the flow to one of the stable equilibria. For
the three player game, the unstable equilibria are the following
three matrices:
$$\left ( \begin{array} {ccc} 0 & {1 \over 2} & {1 \over 2} \\
   0 & 0 & 1 \\ 0 & 1 & 0 \end{array} \right ) \hspace{.5in}
\left ( \begin{array} {ccc} 0 & 0 & 1 \\ {1 \over 2} & 0 & {1 \over 2} \\
   1 & 0 & 0 \end{array} \right ) \hspace{.5in}
\left ( \begin{array} {ccc}  0 & 1 & 0 \\ 1 & 0 & 0 \\
   {1 \over 2} & {1 \over 2} & 0 \end{array} \right ) \; .$$
These correspond to cases where one of the three agents is
entirely ignored, and splits her visits equally between the other
two.  The probability that $p (t)$ is within $\epsilon$ of one of
these traps is roughly $3 \epsilon t^{-1/3}$, so with $t = 1,000$
we find a reasonably high probability that $p (1000)$ is not near
the uniform probability matrix but is instead still near one of
the unstable equilibria.  This persists with reasonable
probability well beyond $t = 10^6$.  For greater population sizes
similar phenomena apply.  Convergence to the invariant set is
relatively slow.  However, for large populations, say 20 or more,
another phenomenon takes place.  The portion of the space of
possible $p$ matrices that are within $\epsilon$ of the possible
limits goes to 1; this is known as the concentration of measure
phenomenon (29).  Thus it becomes very unlikely to get stuck
initially far away from the limit, simply because the initial
randomness will very likely lead to a point very near a possible
limit.  Thus for large populations, the dynamics appear very
similar to the dynamics for Friends I.

\section{Making enemies} \label{ss:negative}

Let us change the ``Making Friends'' model in just one way.
Instead of being rewarded, agents are punished; instead of
uniformly positive interactions, we have uniformly negative ones:

\hspace{.5in}
\begin{tabular}{c|c|}
Enemies I & Host \\
\hline
Visitor & $(-1,0)$ \\
\hline
\end{tabular}
\hspace{.5in}
\begin{tabular}{c|c|}
Enemies II & Host \\
\hline
Visitor & $(-1,-1)$ \\
\hline
\end{tabular}

Instead of interactions being reinforcing, we take them as
inhibiting.  The dynamics of inhibition might be modeled in a
number of ways.  Continuing to use the update
equation~(\ref{eq:update weights}) will not work because the
weights will end up becoming negative and the visitation
probabilities in equation~(\ref{eq:prob1}) will be meaningless.
In this section we explore two other possible rules for updating
the weights so as to inhibit past behavior.  With negative
reinforcement, it is easy to predict what will happen: the social
network always becomes uniform, and the dynamics are not sensitive
to the particular updating mechanism.  Indeed this is what
happens.  Since there are no surprises, and since this model is
just a building block for a model with both structural and
strategic dynamics, we keep the discussion brief.

\subsection{The transfer model}

Consider a three player model with the following update rule on
the weights.  Initial weights are all positive integers.  When $i$
visits $j$, the weight $w_{ij}$ is diminished by 1 and the weight
$w_{ik}$, $k \neq i,j$, is increased by 1.  This is equivalent to
the Ehrenfest model of heat exchange between two bodies (30).  In
the original Ehrenfest model there are two urns.  A ball is drawn
at random from among all balls in both urns and transferred to the
other urn.  The distribution of balls tends to the binomial
distribution, where each ball is independently equally likely to
be in either urn.  In Making Enemies, with transfer dynamics and
three players, each player may be thought of as having such a pair
of urns.  The urns are independent.

Since the number of balls is fixed, an Ehrenfest urn is a Markov
chain with a finite number of states, where the states consist of
distributions over the two urns.  For example, if there are only
two balls, then there are three states, $S1$, $S2$ and $S3$,
corresponding to urn cardinalities of $(2,0)$, $(1,1)$ and
$(0,2)$.  The transition matrix for this Markov chain is
$$\left ( \begin{array} {ccc} 0 & 1 & 0 \\
   {1 \over 2} & 0 & {1 \over 2} \\ 0 & 1 & 0 \end{array} \right )$$
and the unique stationary vector is $(1/4 , 1/2 , 1/4)$.  In
contrast to the P\'olya urn, we do not have convergence of the
conditional probabilities of visits at each stage given the
present: at any time, given the present composition, the
probability of a given visit may be 0 , $1/2$ or 1, depending on
the composition of the urns belonging to the visitor.  However, if
the number of balls, $N$ is large, approximately equal visiting
probabilities are very likely in the following sense.  The
invariant distribution is binomial, which is concentrated around
nearly even distributions when the number of balls is large.  Thus
with high probability, no matter what the initial state, after
roughly $N \log N / 2$ steps (31), the composition of an urn with
$N$ balls will be close to a draw from a binomial distribution.
The conditional probability of either of the two possible visits,
will therefore be close to $1/2$, and will tend to remain there
with high probability. Kac (32) uses these properties to resolve
the apparent paradoxes that beset Bolzmann's discussion of
irreversibility in statistical mechanics.

\subsection{The resistance model}

The transfer model allows for a finite cumulative amount of
negative reinforcement, and indeed yields a finite Markov chain.
Let us explore a rather different model, termed the {\em
resistance model}, in which negative payoffs generate resistance.
Initially every choice has resistance 1. The magnitude of a
negative payoff is added to its associated resistance, so the
equation~(\ref{eq:update weights}) becomes
$$w_{ij} (t+1) = w_{ij} (t) + |u (i,j;t)| \, .$$
In the case at hand, when all payoffs are negative, the
probability of $i$ visiting $j$ is proportional to the reciprocal
of the resistance:
$$p_{ij} = \prob (\mbox{agent } i \mbox{ visits } j) =
   {1 / w_{ij} \over \sum_{k=1}^n 1 / w_{ik}}$$
with $1 / w_{ii} = 0$ by convention.  The dynamics of Enemies I
and Enemies II under resistance dynamics are easy to describe.

\begin{theorem} \label{theorem:enemies}
For Enemies I or Enemies II, from any initial conditions, the
probability matrix $p (t)$ converges to the uniform probability
matrix $\unif$ where $\unif_{ij} = 1 / (n-1)$ for any $i \neq j$.
The of convergence is rapid: of order $N\log N$ if the initial
resistances are of order $N$.  The deviations from uniform obey a
central limit theorem:
$$t^{1/2} (p - \unif) \to X$$
where $X$ is a multivariate normal with covariance matrix of rank
$n(n-1)$ in Enemies I and $n(n-1)/2$ in Enemies II.  In other
words, deviations from uniformity are independent normals, subject
to the constraints of adding up to zero for each individual and,
in the case of Enemies II, the constraints of symmetry.
\end{theorem}

The central limit theorem may be derived from a stronger,
functional central limit theorem, linearizing the system near the
uniform probability to see that the paths
$$t \mapsto N^{-1/2} (p (N t) - \unif)$$
converge in distribution as $N \to \infty$ to a multivariate
Ornstein-Uhlenbeck process.  The rate of convergence follows from
standard coupling arguments.

While uniform positive reinforcement breeds structure from
unstructured initial conditions, uniform negative reinforcement
evidently breeds uniformity even from structured initial
conditions.  It would appear, therefore, that the customary random
encounter (mean-field) model is more suitable for Making Enemies
than Making Friends.

\subsection{A better model?}

We would like a model that allows for both positive and negative
reinforcement.  A natural choice is to let $w_{ij}$ keep track of
the log-likelihood for $i$ to visit $j$, so that probability of
$i$ visiting $j$ is given by
\begin{equation} \label{eq:log update}
p_{ij} = \prob (\mbox{agent } i \mbox{ visits } j) =
   {\exp (w_{ij}) \over \sum_{k=1}^n \exp (w_{ik})} \, .
\end{equation}
In the next section we will see a property this rule has in common
with rules that discount the past, namely that it leads to being
trapped in a deterministic state where $i$ always visits the same
$j$.

\begin{question}
Is there a model incorporating both positive and negative
reinforcement, that is realistic, tractable, and non-trapping?
\end{question}

\section{Perturbations of the models}

In this section we add two features, noise and discounting,
commonly used to create more realistic models.  We examine the
effects on social structure.  In particular, these lead to varying
degrees of subgroup formation.

\subsection{Discounting the past}

In the foregoing models, a positive (or negative) payoff in the
distant past contributes equally to the weight (or resistance)
assigned to an edge as does a like payoff in the immediate past.
This is implausible, both psychologically and methodologically. As
a matter of psychology, memories fade.  From the standpoint of
inductive logic, it is not at all certain that the learner is
dealing with stationary probabilities - indeed, in cases of prime
interest she is not.  For this reason, recent experience may have
a better chance of being a relevant guide to future action than
the remote past.

A simple and standard way to modify the models to reflect this
concern is to introduce discounting of the past.  We will
concentrate here on the models of Making Friends.  After each
interaction we will now multiply the weights of the previous stage
by a discount factor, $d$, between 0 and 1.  The we add the
undiscounted payoffs from the present interaction to get new
weights.  The modification of the dynamics has a dramatic effect
on the Making Friends models.

For Friends I, it is immediately evident from simulations with $d
= .9$, say, and ten players, that the probabilities $p_{ij}$
converge to 0 or 1. In other words, each individual ends up always
visiting the same other individual.

In Friends II, simulations show the group breaking into pairs,
with each member of a pair always visiting his or her ``partner''.
Which pairs form depends on the randomness in the early rounds of
visits, but pairs always form.  In fact there are other possible
limit states, but their frequency is low except at more extreme
discount rates.  The set of possible limit states may be described
as follows.  Some agents are grouped in pairs, each member of a
pair always visiting the other.  Other agents are grouped in {\em
stars}.  These are clusters of size at least three, in which one
agent, called the {\em center}, visits each of the others with
positive frequency, while the others always visit the center.

\subsection{Analysis of discounting the past}

It is worth giving a rigorous derivation of the above behavior,
since it will shed some light on a defect in the most obvious
log-likelihood model to incorporate positive and negative
reinforcement.  Our derivation highlights this, although the
results for discounted Friends I may also be derived from a
theorem of H. Rubin (33, Page 227).

\begin{theorem} \label{theorem:discounted friends}
In Friends II with discount rate $d < 1$, there is always a
partition into pairs and stars and a random time after which each
member of a pair visits only the other member of the pair and each
non-central member of a star visits only the center. In Friends I,
there is a random function $f$ and a random time after which each
player $i$ always visits $f(i)$.
\end{theorem}

\noindent{\em Sketch of Proof:} The analysis for Friends I is
similar but easier, so we prove the statement only for Friends II.
With each probability matrix $p$ we associate a graph $G(p)$ as
follows.  The edge $(i,j)$ is in the graph $G$ if the probability
$p_{ij} > \epsilon$, where $\epsilon < 1/(2n)$ is some fixed positive
number. Among those graphs having at least one edge incident to
each vertex, let $S$ denote the minimal such graphs, that is, ones
for which deleting any edge results in an isolated vertex.  It is
easy to see that $S$ is the set $G(p)$ for all $p$ satisfying the
conclusion of the theorem.

The principle behind the analysis of discounted Friends is that
the future behavior of $p$ is largely determined by the present
$G(p)$.  In particular, we find a $\delta > 0$ such that from
any state $p$, for each subgraph $H$ of $G(p)$ such that $H \in
S$, there is a probability at least $\delta^2$ that for all
sufficiently large $t$, $G(p(t)) = H$.  We show this in two steps:
(1) with probability at least $\delta$, there is some $t$ for
which $G(p(t)) = H$; (2) from any state $p$ such that $G(p) = H$,
there is probability at least $\delta$ that $G(p(t))$ is equal to
$H$ for all later times, $t$.

To see why~(1) is true, for $H \in S$, let $f_H$ be any function
on vertices of $H$ for which each value $f(i)$ is a neighbor if
$i$. Observe that there is a number $k$ such that from any state
$p$ with $H \subseteq G(p)$, if each vertex $i$ visits $f(i)$ for
the next $k$ rounds, then $G(p(k)) = H$.  For each round of
visits, this probability is at least $\epsilon^n$, where $n$ is
the number of vertices, so taking $\delta \leq \epsilon^{kn}$
establishes~(1).  For~(2), it suffices to show that with
probability $\delta$ each agent visits a neighbor in $H$ at all
later times.  For each agent $i$, the sum over $j$ not neighboring
$i$ in $H$ of $p_{ij}$ is at most $n \epsilon < 1/2$ by the
definition of $G(p) = H$.  After $k$ rounds of visits where agents
only visit their neighbors in $H$, this must decrease to at most
$(1/2) d^k$.  Thus the probability of $N$ rounds of visits only to
neighbors in $H$ is at least
$$\prod_{k=0}^{N-1} \bigg ( 1 - {1 \over 2} d^k \bigg )^n .$$
Sending $N$ to infinity yields a convergent infinite product,
since $(1/2) d^k$ is summable.  Taking $\delta$ to be less than
the infinite product proves~(2).

With~(1) and~(2), the rest is a standard tail argument.  The
constraints on evolution are such that $G(p(t))$ always contains
at least one graph in $S$.  As long as it contains more than one
graph in $S$, there is always a probability of at least $\delta$
of permanently settling into each one.  Thus, with probability 1,
eventually $G(p(t))$ is equal to some $H \in S$ for all future
times.  This is equivalent to the conclusion of the theorem.
$\qed$

\noindent{\em Remark:} It is actually shown that in~(2), if we
choose $\epsilon$ sufficiently small, we can choose $\delta$
arbitrarily close to 1.

We now also see why the log-likelihood rule~(\ref{eq:log update})
leads to fixation of a degenerate structure.  Under these
dynamics, an equivalent phenomenon occurs to~(1) in the proof of
Theorem~\ref{theorem:discounted friends}.  For a pair $(i,j)$
whose interaction has a positive mean, if the pair plays
repeatedly, we will see $w_{ij} (t) / t \to \mu > 0$.  The
probability the $i$ will ever switch partners, once having tried
$j$ a few times is at most on the order of
$\sum_{k=0}^\infty B \exp (- k \mu)$, where
$B = \exp (\sum_{l \neq j} w_{il})$.  From here it is easy to
construct an argument parallel to the proof of
Theorem~\ref{theorem:discounted friends}, to show that in presence
of a game with positive mean payoff, discounted structural
dynamics lead with probability 1 to fixation at a pairing.

\subsection{Introduction of noise}

A common feature in models of adaptation is the introduction of
noise: a small chance of a behavior other than the one chosen by
the dynamical equation for the model.  This may stem from an
agent's uncertainty, from agent error, or from circumstances
beyond an agent's control.  Alternatively, an agent may
purposefully add noise to her strategy in order to avoid becoming
wedded to a less than optimally efficient strategy or structure.

 From a methodological point of view, noise that does not go to
zero with time transforms the model into an ergodic Markov chain.
No state is then trapping.  To the extent that the trapping states
produced by discounting or linear log-likelihood are unrealistic,
we may hope to mitigate the problem by adding a noise component.
Since dynamics with a noise term do not lead to a single state,
the outcome is usually phrased in terms of {\em stochastically
stable states} (34).  A state is termed stochastically stable if
the chance of finding the system near that state does not go to
zero as the magnitude of the noise term goes to zero.

Neither discounting nor noise will affect the limiting behavior of
Making Enemies.  For Making Friends, let us modify the probability
rule~(\ref{eq:prob1}) so that in the $n$-player game, the
probability of $i$ visiting $j$ is now some fixed positive number
$\epsilon / (n-1)$, plus $(1-\epsilon)$ times what it was before:
$$p_{ij} = {\epsilon \over n-1} + (1 - \epsilon) {w_{ij}
   \over \sum_k w_{ik}} \, .$$
The effect of this is to push the system by $\epsilon$ toward the
uniform point $\unif$.  Neither Friends I nor Friends II is now a
martingale, and the stable set of each is reduced to the single
point $\unif$.  Since this is true at any noise level $\epsilon >
0$, we see that there is only one asymptotically stable point.
Since the qualitative outcome is sensitive to the existence of a
noise term, it is incumbent to ask with regard to specific models
whether a noise term is natural and realistic.

\subsection{Noise and discounting}

In the presence of a discount $d < 1$ and a noise term $\epsilon >
0$, if $1-d$ is much smaller than $\epsilon$ then the discount is
so low that the noise term wipes out any effect the discounting
might have had.  In the other case, where $d$ is held fixed and
$\epsilon$ tends to zero, we may ask about the asymptotically
stable states of system with past discounting dynamics.  For
Friends I, nothing much interesting happens: discounting causes
the limiting state to be degenerate; with noise, the system may
jump from one such state to the other, which does not change which
states are stochastically stable.

For Friends II, as long as the number of players $n$ is at least
4, the introduction of noise does indeed change the set of
stochastically stable states: it gets rid of stars.  Simulations
show that pairings are by far the most prevalent states in
discounted Friends II, with a star of size 3 forming when
necessitated by an odd number of players.  We now show that states
with more than one star, or a star of size greater than 3, are not
stochastically stable.

\begin{theorem} \label{theorem:stable}
In Friends II, with discounting, with $n$ players, and with noise
tending to zero, the stochastically stable states are those that
are either unions of pairs (if $n$ is even) or pairs plus a single
star of size 3 (if $n$ is odd).
\end{theorem}

\noindent{\em Sketch of Proof:}  Let $S$ denote the graphs
corresponding to possible limit states as in the proof of
Theorem~\ref{theorem:discounted friends}, and let $S_0 \subseteq
S$ denote those graphs with no stars (perfect pairings) or with a
single star of size 3.  The important features of the relation of
$S$ to $S_0$ are as follows.  (1) if $G$ is the result of adding a
single edge to a graph in $S_0$, then $G$ contains no graph in $S
\setminus S_0$.  (2) for any $G \in S$ there is a chain $G = G_1 ,
G_2 , \ldots , G_k$ leading to $S_0$, where each $G_{j+1}$ may be
obtained from $G_j$ be adding an edge and then deleting two edges.
Property~(1) is apparent.  To verify~(2), note that if $H \in S$
and $i$ and $j$ are non-central vertices in stars of $H$, and they
are not both in the same star of size 3, then adding the edge
between $i$ and $j$ and removing the two edges previously incident
to $i$ and $j$ produces a new graph in $S$.  Iterating this
procedure starting from $H = G_1$ leads in finite time (since the
number of edges decreases each time) to an element of $S_0$.

We now follow the usual method for determining stochastic
stability (35).  Let the probability $\rho$ of disobeying the
structural dynamics equation~(\ref{eq:prob1}) be very small. If
$\epsilon$ (in the definition of $S$) is very small, then a state
$p$ with $G(p) = G \in S$ will have $G(p(t)) = G$ for all later
times with high probability, until there is a disobeying move.
After a single disobedience, the graph $G(p)$ will be the union of
$G$ with one extra edge.  By the remark after the proof of
Theorem~\ref{theorem:discounted friends}, we see that after a
disobedience, the graph will then relax to some subgraph in $S$.
By property~(1), if $G \in S_0$ then this subgraph is again in
$S_0$.  Thus a single disobedience followed by relaxation back to
$S$ will never escape $S_0$.  Hence the probability of jumping to
$S \setminus S_0$ is of order $ \rho^2$, which implies that states
in $S_0$ stay in $S_0$ for time at least $\rho^{-2}$. On the other
hand, by property~(2), from any state in $S \setminus S_0$, there
is a chain of single disobediences, such that allowing the system
to relax after each may with positive probability land you back in
$S_0$.  Thus the expected time spent in $S \setminus S_0$ before
returning to $S_0$ is at most of order $\rho$.  Thus the process
spends $(1 - \rho)$ portion of the time in $S_0$, and sending
$\rho$ to zero, we see that only states in
$S_0$ are stochastically stable.  It is easy to see that all of
these are indeed stochastically stable.  $\qed$

\section{Reinforcement by games of nontrivial strategy}

So far we have only considered a baseline model of uniform
reinforcement, which turned out still to have nontrivial
structural behavior.  Now we examine a reinforcement scheme
resulting from the payoff of a nontrivial game.  We will consider
the case where evolution of strategy is slower than evolution of
structure.  Thus, we will consider the agents as divided into
types, each type always playing a fixed strategy, and see what
sort of interaction structure emerges.  We then extend this by
allowing strategic switching of types.  We find that coordination
of strategy occurs, though whether players coordinate on the
risk-dominant or payoff-dominant strategy depends on parameters of
the model such as the rate of strategic evolution.  Depending on
conditions of the model, the social network may or may not split
up into pairs.

\subsection{Rousseau's Stag Hunt}

Consider a two-player version of Rousseau's {\em Stag Hunt} (36).
The choices are either to hunt stag or to hunt rabbit (hare, in
the original).  It takes two person cooperating to effectively
hunt a stag, while one person acting independently can hunt a
rabbit. Bagging a stag brings a greater payoff.

\begin{tabular}{c|c|c|}
& Hunt Stag & Hunt Rabbit \\
\hline
Hunt Stag & $(1,1)$ & $(0,.75)$ \\
\hline
Hunt Rabbit & $(.75 , 0)$ & $(.75 , .75)$ \\
\hline
\end{tabular}

There are two equilibria in this game: both hunt stag, and both
hunt rabbit.  The first carries the higher payoff and is said to
be {\em payoff dominant}; the second carries the least risk and is
said to be {\em risk dominant} (37).  In models without structural
dynamics, Kandori, Mailath and Rob (38) have shown that only the
risk dominant equilibrium of a two player coordination game is
stochastically stable.  In the presence of structural dynamics, we
will describe a more optimistic conclusion.

\begin{theorem} \label{theorem:stags}
Suppose Stag Hunt is played by $2n$ players, with structural
dynamics given by equation~(\ref{eq:prob1}) and cumulative
weighting dynamics~(\ref{eq:update weights}) with no noise or
discounting. Then in the limit, stag hunters always visit stag
hunters and rabbit hunters visit rabbit hunters.
\end{theorem}

\noindent{\em Sketch of Proof:} First note that no visit of a stag hunter
to a rabbit hunter is ever reinforced.  Thus $w_{ij} (t) = 1$ for
all $t$ if $i$ is a stag hunter and $j$ is a rabbit hunter.
Observing that the weights $w_{ij} (t)$ go to infinity when $i$
and $j$ are both stag hunters, we see that the probability of a
stag hunter visiting a rabbit hunter goes to zero.

Next consider the subpopulation of rabbit hunters, call it $A$.
For $i \in A$, let
$$Z(i,t) = {\sum_{j \notin A} w_{ij} \over \sum_{j=1}^n w_{ij}}$$
denote the probability of visiting a given rabbit hunter visiting
a stag hunter on the next turn.  The expected value of $Z(i,t+1)$
changes according to the formula
$$\expect (Z(i,t+1) | Z(i,t)) = Z(i,t) + t^{-1} Y(i,t)$$
where $Y(i,t)$ is the proportion of increase in expected weight
$w_{ij}$ due to $j \notin A$:
$$Y(i,t) = {\sum_{j \notin A} p_{ij} + p_{ji} \over \sum_{j=1}^n
   p_{ij} + p_{ji}} \, .$$
Ignoring the terms $p_{ji}$ in both the numerator and denominator
of the above expression would lead to exactly $Z(i,t)$.  The terms
$p_{ji}$ for $j \notin A$ are known to be small, while the total
from the terms $p_{ji}$ for $j \in A$ cannot be small.
Consequently, $Y(i,t) < (1 - \epsilon) Z(i,t)$ for some $\epsilon
> 0$, whence
$$\expect (Z(i,t+1) - Z(i,t) | Z(i,t)) \leq
   - {\epsilon Z(i,t) \over t} \, .$$
Since the increments in $Z(i,t)$ are bounded by $C/t$, there are a
$\lambda , \mu > 0$ for which $\exp (\lambda Z(i,t) + \mu \log t)$ is a
supermartingale, which implies that $Z(i,t)$ converges to zero
exponentially fast in $\log t$.  $\qed$

Introduction of a discount rate changes this outcome. Stag hunters
still end up visiting stag hunters, since even discounted
reinforcement beats a reinforcement of zero, but now rabbit
hunters will get locked either into pairs and stars as in Making
Friends, or into repeated visits to a single stag hunter. These
limit states are all invariant under introduction of noise.  When
a rabbit hunter visits a stag hunter the loss to society is the
$.75$ that another rabbit hunter would have profited from the
visit.  The model is evidently weak here, since it allows only one
visit {\em by} each agent but any number of visits {\em to} each
agent in a round of visits.  That is, a more realistic loss would
be the stag hunter's wasted time when visited by the rabbit
hunter.

It should be noted that although the stochastically stable states
include ones that are not optimally efficient, the optimally
efficient states (those states where rabbit hunters visit rabbit
hunters) will have an edge.  Due to the possibility of reciprocal
reinforcement, it will be easier for a rabbit hunter to switch
from visiting a stag hunter to visiting a rabbit hunter, than {\em
vice versa}. Secondly, when the discount rate is near 1, the model
behaves like the undiscounted model for a long enough time that it
is very unlikely for a rabbit hunter to get locked into visiting a
stag hunter in the first place. Simulations of Stag Hunting with
ten players and $d = .9$, seem to show that rabbit hunters
``always'' visit rabbit hunters.  Due to both of the effects
mentioned above, the system is nearly always found in an optimally
efficient state, even though there are stochastically stable
states that are not optimally efficient.

\subsection{Co-evolution of structure and strategy}

To the previous model, we now add the possibility of an agent
switching states: a stag hunter may decide to become a rabbit
hunter, or a rabbit hunter may become bold and hunt stag.  When
this kind of strategic evolution is faster than the structural
evolution, we know from studies of random encounter models that
the risk dominant equilibrium of everyone hunting rabbits will be
arrived at while the network is still near its initial state of
uniform visitation probabilities.

Whether strategic dynamics are faster or slower than structural
dynamics depends, of course, on the activity being modeled;
sometimes interaction structure is externally imposed, while
sometimes it is more easily modified than strategy or character.
Let us suppose that the investment in re-training as a different
kind of hunter is great, so between each round of visits there is
only a small chance that one of the hunters will change types.
Then we have seen that hunters always (with no noise or
discounting) or nearly always (in discounted models) hunt with
others of like type.  This eliminates the risk inherent in random
encounters, and allows hunters to profit from switching to stag
hunting after an initial period where they find another stag
hunter.  Slow strategic adaptation gradually converts rabbit
hunters to stag hunters and the payoff dominant strategy
dominates.

We describe here the results of simulations of Stag Hunting for
1,000 time steps, where with some probability $q$ at any given
time, an individual changes type to whichever type was most
successful in the previous round.  When $q = .1$, we found that in
22\% of the cases all hunters ended up hunting stag, while in 78\%
of the cases, all hunters hunted rabbit.  Thus there was perfect
coordination, but usually not to the most efficient equilibrium.
On the other hand, when $q = .01$, the majority (71\%) of the
cases ended in the optimal state of all hunting stag, while 29\%
ended up all hunting rabbit.  Increasing the initial edge weights
made it far less likely to reach the stag hunting equilibrium,
since stag hunters took a long time to perfectly align, and
without alignment, the previous round's best strategy was almost
always rabbit hunting.  For instance, if the initial weights were
1,000 for each visit, under 1\% of the cases ended up all stag
hunting, whether $q$ was $.1$ or $.01$.

Once hunters largely cease to visit hunters of opposite type, the
structural evolution within each of the two subpopulations is a
version of Friends II.  The resulting social structure will not be
a perfect pairing, but will have each rabbit (stag) hunter
visiting each other rabbit (stag) hunter, but with varying
probabilities.

\section{Conclusion}

We have taken some basic steps in exploring dynamics of evolution
of interaction structures and co-evolution of structure and
strategy.  The ultimate goals are to create models that are more
true to life, and to find theoretical bases for observed behaviors
of systems, including prediction of selection between multiple
equilibria.

The particular dynamics we use here are only examples, but it
turns out that the simplest of these may deliver interesting and
surprising results. Even in baseline models where the game being
played is degenerate, we find spontaneous emergence of structure
from uniformity and spontaneous emergence of uniformity from
structure.  We find processes with extremely long transient modes,
where limiting behavior is not a good guide for predicting
behavior after thousands of trials.

The social interaction structures that emerge tend to separate the
population into small interaction groups within which there is
coordination of strategy.  This separation may be complete, as in
discounted Friends II, or may be only a tendency, as in the
non-discounted versions of Friends and Stag Hunting.

When we combine structure and strategy dynamics for a non-trivial
game, the Stag Hunt, we find that the probable outcomes depend on
the timing. Where structure is frozen in a random encounter
configuration we get the expected risk-dominant equilibrium
outcome. But when structure is fluid relative to strategy,
structural adaptation neutralizes the risk and we get the socially
efficient payoff dominant equilibrium. Varying between these
extremes can give one or the other result with different
probabilities - or may leave the group in a state where both
strategies are used.  We expect to see structure dynamics making a
difference in other games as well. Indeed, we have some
preliminary simulation evidence showing this to be true for a
bargaining game (``split the dollar''), and for a simple
coordination game.

There are many more avenues to pursue.  As mentioned in
Section~\ref{ss:negative}, it would be desirable to find a model
in which positive and negative reinforcement are present, but
trapping does not occur. We have not modeled any interaction among
three or more players. We also have yet to model any explicit
interaction between strategy and structure: the choice of a
partner to play with and a strategy to play against that partner
need not be independent.

One could continue adding complexity so as to allow information to
affect structural evolution, to include communication between
players, and so forth.  Our main point is this.  Structural change
is a common feature of the real world. A theory of strategic
interaction must take account of it. There is a mathematically
rich theory which develops relevant tools. We believe that
explicit modeling of structural dynamics, and the interaction of
structure and strategy, will generate new insights for the theory
of  adaptive behavior.

\noindent{We} wish to thank Persi Diaconis,
Joel Sobel and Glenn Ellison, for bringing us together, for
helpful discussions, and for greatly improving our awareness of
the relevant literature.

\setlength{\parindent}{0pt}

1.  Hofbauer, J. and Sigmund, K. (1988).  {\em The theory of
evolution and dynamical systems}  (Cambridge University Press,
Cambridge).

2.  Weibull, J. (1997).  {\em Evolutionary game theory} (MIT
Press: Cambridge, MA).

3.  Feldman, M. and Thomas, E. (1987). Behavior-dependent contexts
for repeated player of the Prisoner's Dilemma II: Dynamical
aspects of the evolution of cooperation. {\em J. Theor. Biol.}
{\bf 128},~297--315.

4.  Wright, S. (1921).  Systems of mating III: assortative mating
based on somatic resemblance.  {\em Genetics} {\bf 6}, 144--161.

5.  Wright, S. (1945).  Tempo and mode in evolution: a critical
review.  {\em Ecology} {\bf 26}, 415--419.

6.  Jackson, M. and Watts, A. (1999).  On the formation of
interaction networks in social coordination games.  {\em Working
paper.}

7.  Durrett, R. and Neuhauser, C. (1997).  Coexistence results for
some competition models.  {\em Ann. Appl. Prob.} {\bf 7}, 10--45.

8.  Kang, H.-C., Krone, S. and Neuhauser, C. (1995).  Stepping
stone models with extinction and recolonization.  {\em Ann. Appl.
Prob.} {\bf 5}, 1025--1060.

9.  Durrett, R. and Levin, S. (1994).  The importance of being
discrete (and spatial).  {\em Theor. Pop. Biol.} {\bf
46},~363--394.

10.  Blume, L. (1993).  The statistical mechanics of strategic
interaction.  {\em Games Econ. Behav.} {\bf 5}, 387--423.

11.  Pollack, G. B. (1989).  Evolutionary stability on a viscous
lattice.  {\em Social Networks} {\bf 11} 175--212.

12.  Lindgren, K. and Nordahl, M. (1994).  Evolutionary dynamics
and spatial games.  {\em Physica D} {\bf 75}, 292--309.

13.  Anderlini, L. and Ianni, A. (1997).  Learning on a torus.
In: {\em The dynamics of norms}, ed. C. Bicchieri, R. Jeffrey, and
B. Skyrms (Cambridge University Press, Cambridge) 87--107.

14.  Ellison, G. (1993).  Learning, local interaction, and
coordination.  {\em Econometrica} {\bf 61}, 1047--1071.

15.  Watts, D. and Strogatz, S. (1998) Collective dynamics of
``small-world'' networks. {\em Nature} {\bf 393},~440--442.

16.  Schelling, T. (1969)  Models of Segregation.  {\em American
Economic Review, Papers and Proceedings} {\bf 59}, 488--493.

17.  Schelling, T. (1971).  Dynamic models of Segregation. {\em
Journal of Mathematical Sociology} {\bf 1}, 143--86.

18.  Epstein, J. and Axtell, R. (1996). {\em Growing Artificial
Societies}  (MIT/Brookings: Cambridge, MA).

19.  Zajonc, R.B. (1968). Attitudinal effects of mere exposure.
{\em Journal of Personality and Social Psychology Monograph} {\bf
9}, 1--28.

20.  Coppersmith, D. and Diaconis, P. (1987).  Reinforced Random
Walk. {\em Unpublished manuscript.}

21.  Arthur, W. B. (1989).  Competing technologies, increasing
returns, and lock-in by historical events.  {\em Econ. J.} {\bf
99},~116--131.

22.  Eggenberger, F. and P\'olya, G. (1923).  \"Uber die Statistik
verketter Vorgange.  {\em Zeit. Angew. Math. Mech.} {\bf 3},
279--289.

23.  Johnson, N. and Kotz, S. (1977).  {\em Urn models and their
application} (John Wiley \& Sons, New York).

24.  Durrett, R. (1996).  {\em Probability: Theory and examples,
2nd edition} (Duxbury Press, Wadsworth Publishing Company,
Belmont, CA).

25.  Pemantle, R. and Skyrms, B. (2000).  Reinforcement schemes
may take a long time to exhibit limiting behavior.  {\em In
Preparation.}

26.  Robbins, H. and Monro, S. (1951).  A stochastic approximation
method.  {\em Ann. Math. Statist.} {\bf 22}, 400--407.

27.  Pemantle, R. (1990).  Nonconvergence to unstable points in
urn models and stochastic approximations.  {\em Ann. Probab.} {\bf
18}, 698--712.

28.  Bena\"im, M. and Hirsch, M. (1995).  Dynamics of Morse-Smale
urn processes.  {\em Ergodic Theory and Dynamical Systems} {\bf
15}, 1005--1030.

29.  Talagrand, M. (1995).  Concentration of measure and
isoperimetric inequalities in product spaces.  {\em IHES Publ.
Math.} {\bf 81}, 73--205.

30.  Ehrenfest, P. and Ehrenfest, T. (1907).  \"Uber zwei bekannte
Einwande gegen das Bolzmannshe H-Theorem.  {\em Phys. Zeit.} {\bf
8}, 311--314.

31.  Diaconis, P. and Stroock, D. (1991).  Geometric bounds for
eigenvalues of Markov chains.  {\em Ann. Appl. Prob.} {\bf 1},
39--61.

32.  Kac, M. (1947).  Random walk and the theory of Brownian
motion. {\em Amer. Math. Monthly} {\bf 54}, 369--391.

33.  Davis, B. (1990).  Reinforced random walk.  {\em Prob. Th.
Rel. Fields} {\bf 84}, 203--229.

34.  Foster, D. and Young, H. P. (1990).  Stochastic evolutionary
game theory.  {\em Theor. Pop. Biol.} {\bf 38}, 219--232.

35.  Ellison, G. (2000).  Basins of attraction, long run
stochastic stability, and the speed of step-by-step evolution.
{\em Rev. Econ. Studies, to appear.}

36.  Rousseau, J.-J. (1984).  {\em A discourse on inequality}.
Tr. Maurice Cranston (Penguin Books, London).

37.  Harsanyi, J. and Selten, R. (1988).  {\em A general theory of
equilibrium in games}  (MIT Press, Cambridge, MA).

38.  Kandori, M., Mailath, G. and Rob, R. (1993).  Learning,
mutation, and long run equilibria in games.  {\em Econometrica}
{\bf 61}, 29--56.

\end{document}